\documentclass[12pt]{amsart}
\usepackage{geometry}                
\usepackage{graphicx}
\usepackage{amssymb}
\usepackage{eqnarray,amsmath}

\usepackage{epstopdf}
\usepackage{enumerate}

\theoremstyle{plain} 
\newtheorem{theorem}{Theorem}

\theoremstyle{definition} 
 
\newtheorem{conjecture}{Conjecture} 
\newtheorem{claim}{Claim}

\newtheorem{problem}{Problem}

\theoremstyle{remark}

\title{The  (7,4)-conjecture in finite groups  }

\author[1]{J\'ozsef Solymosi}

\address{Department of Mathematics, University of British Columbia \\ 1984 Mathematics Road, Vancouver, BC, V6T 1Z4, Canada}
\email{solymosi@math.ubc.ca}
\thanks{Research was supported by NSERC, ERC-AdG. 321104, and OTKA NK 104183 grants.}
\begin{document}

\begin{abstract} The first open case  of  the Brown, Erd\H{o}s, S\'os conjecture is equivalent to the following; For every $c>0$ there is a threshold $n_0$ so that if a quasigroup has order $n\geq n_0$ then for every subset of triples of the form $(a,b,ab),$ denoted by $S,$ if $|S|\geq cn^2$ then there is a seven-element subset of the quasigroup which spans at least four triples of the  selected subset  $S.$  In this paper we prove the conjecture for finite groups. 
\end{abstract}

\maketitle

\section{Introduction}
This paper is about  proving a special case  of a famous conjecture in extremal combinatorics. The conjecture originates from  Brown, Erd\H{o}s, and T. S\'os \cite{BES}. Before we state it let us introduce some notations we are going to use. {\em Triple systems}  are families of three-element subsets of a finite set.  In the theory of hypergraphs such systems are called three-uniform hypergraphs.  If a triple system has many triples, if it is dense in some sense, that is a {\em global} property. Usually it is hard to show that  dense systems have {\em locally} dense subsystems. For example Tur\'an's conjecture states that if the number of triples is more than $\frac{5}{9}{n\choose 3}$ in a triple system $\mathcal{T}$  on $n$ elements then there are four elements that all four triples spanned by the four elements are in $\mathcal{T}.$ (The 3-uniform hypergraph contains a clique, $K_4^{(3)}.$) A more general question is the following; What can we say about the density of a triple system if one knows that no $k$ elements span $\ell$ or more triples?  Depending on the values of $\ell$ and $k$ the question might be a very hard one. Understanding how global properties induce local properties is a central problem in combinatorics. 

\medskip
The Brown-Erd\H os-S\'os  conjecture is that for any fixed $k\geq 3$ all triple systems on $n$ elements in which no $k+3$ elements span $k$ triples should be sparse, i.e. it has $o(n^2)$ triples. Note that here sparseness is relative to the fact that such systems have $O(n^2)$ triples. Observe that if $k$ triples have a common 2-element intersection than $k+2$ elements span $k$ triples. Therefore, if no $k+2$ elements carry $k$ triples then the number of triples is at most $3(k-1){n\choose 2}.$  In this paper we will suppose that our triple system is a linear hypergraph, that is no triples share more than one element. If no $k+3$ elements carry $k$ triples then a constant fraction of the triples form a linear hypergraph. 
\begin{claim}
To prove or disprove the Brown-Erd\H os-S\'os  conjecture it is enough to check it for linear 3-uniform hypergraphs.
\end{claim}
To see the claim let's order the $m$ triples and check them one by one. Take the first triple, then select the second if it has at most one common element with the first one, and select the third one if it has at most one common element with the previously selected triples, so on. Every selected triple had at most $k$ triples with 2 elements common, so at the end of the selection we still have at least $m/(k+1)$ triples selected so that no two have two common elements.     \qed

\medskip
In the other direction it was noted in \cite{BES} that a random construction shows that for every $k\geq 3$ there is a $c_k> 0$ such that one can find triple systems with $c_kn^2$ triples on $n$ elements that no $k+2$ elements  span $k$ triples. ($n$ can be chosen arbitrary large) For the sake of completeness we sketch the random construction here. The details can be found in \cite{BES}.

\medskip
\noindent
{\bf Construction:} Choose triples out of the possible $n\choose 3$ triples in an $n$ element set independently, at random, with probability $\delta n^{-1}.$ If in this triple system for some $k+2$ elements there are $k$ or more triples spanned by the elements then remove all such spanned triples from the system. There is a constant $c_k>0$ such that for any choice of $k+2$ elements, the probability that we selected at least $k$ triples out of the possible ${k+2}\choose 3$  is less then  $c_k\delta^kn^{-k}.$ By the linearity of expectations the expected number of the removed triples is less than 

$$c_k\delta^kn^{-k}{n\choose {k+2}}{{k+2}\choose 3}\leq c_k'\delta^kn^2, $$
for some $c_k' > 0$ which depends on $k$ only. If we choose $\delta$ small enough that 

$$c_k'\delta^k\leq \frac{\delta}{2}{n\choose 3}$$ then less than half of the selected triples were removed, so there are still some $c_k''n^2$ triples remain. 

\medskip

One might think that the $k+2,k$ case is solved since the triple systems without $k$ elements carrying $k+2$ triples can not have more than $C_kn^2$ triples on $n$ elements and as the previous construction shows there are such systems with $c_kn^2$ triples. But there is an interesting question which still remains open; The two constants are far apart. In the previous arguments $C_k\rightarrow \infty$ and $c_k\rightarrow 0$ as $k\rightarrow \infty.$ 

\begin{problem} Is it true that for every integer, $k\geq 100$,  if a triple system on $n\geq k+2$ elements contains at least $n^2/100$ triples then it contains $k+2$ points carrying at least $k$ triples? (Of course 100 is just an arbitrary number here. Does the statement hold for some constant?)
\end{problem}

I heard the problem above from Nati Linial first, but probably others had similar questions too. In a related conjecture of Erd\H os - which would imply a negative answer to the previous problem - the question is formulated as follows.

\begin{conjecture}[Erd\H os' Steiner Triple System Conjecture]
For every $r\geq 4$ there are arbitrary large Steiner triple Systems where no $r+2$ elements carry at least $r$ triples.
\end{conjecture}

There are partial results on Erd\H os' Steiner Triple System Conjecture. We refer to the papers \cite{FGG} and \cite{FU} for further details.

\section{Main result}
We will reformulate the Brown-Erd\H os-S\'os  conjecture as a statement in quasigroups. Our hope is that some tools from algebra can be used to attack this notoriously  hard problem.

It was a conjecture of Lindner  \cite{Li} that any partial Steiner triple system of order $u$ can be embedded in a Steiner triple system of order $2u$. This was proved by  Bryant and Horsley in \cite{BH}. 
Any Steiner triple system defines a quasigroup. So, the following is equivalent to the original Brown-Erd\H os-S\'os  conjecture.

\begin{conjecture}[Brown, Erd\H{o}s, S\'os]\label{BESconj}
For every $c>0$ there is a threshold $n_0$ so that if a quasigroup has order $n\geq n_0$ then for every subset of triples of the form $(a,b,ab),$ denoted by $S,$ if $|S|\geq cn^2$ then there is a seven-element subset of the quasigroup which spans at least four triples of the  selected subset  $S.$ 
\end{conjecture}

Now the question is that for which families of quasigroups can we prove the conjecture. The main result of this paper is to show that the (7,4)-conjecture holds for finite groups. 

\begin{theorem}[The Brown-Erd\H{o}s-S\'os Conjecture for groups]\label{main}
For every $c>0$ there is a threshold $n_0$ so that if a group has order $n\geq n_0$ then for every subset of triples of the form $(a,b,ab),$ denoted by $S,$ if $|S|\geq cn^2$ then there is a seven-element subset of the group which spans at least four triples of the  selected subset  $S.$ 
\end{theorem}

In addition to the algebraic techniques there are some combinatorial tools which can be used when one works with triple systems. The most powerful one is the so called Hypergraph Removal Lemma \cite{G,NRS}, which we are going to apply here. The simplest case, the Triangle Removal Lemma,  states that for every dense subset of triples of the form $(a,b,ab)$  there is a six-element subset of the quasigroup which spans at least three triples from the  selected subset.   This is called the (6,3)-theorem. It was proved by Ruzsa and Szemer\'edi \cite{RSz}. In search for the proof of the (7,4) conjecture, Frankl and R\"odl proved the Removal Lemma for 3-uniform Hypergraphs \cite{FR}. We are going to use the following form of the result;

\begin{theorem}[Frankl-R\"odl]\label{FR1}
Let $H_n^{(3)}$ be a 3-uniform hypergraph on $n$ vertices with the property that every edge is the edge of exactly one clique, $K_4^{(3)}.$ Then the number of edges is $o(n^3).$ (For every $\varepsilon >0$ there is a threshold $n_0=n_0(\varepsilon)$ so that if $n\geq n_0$ and $H_n^{(3)}$ has the above property then it has at most $\varepsilon n^3$ edges) 
\end{theorem}

Our application of the above theorem is similar to the technique we used in \cite{So1}. Theorem \ref{FR1} is enough to prove the (7,4) conjecture in groups, however for some quantitative results the following stronger statement is useful.

\begin{theorem}[Frankl-R\"odl]\label{FR2}
For every real number $c>0$ there is a $c'>0$ such that the following holds. If $H_n^{(3)}$ is a 3-uniform hypergraph on $n$ vertices with the property that it has at least $cn^3$ edge-disjoint $K_4^{(3)}$ cliques then  it contains at least $c'n^4$ distinct (but not necessary edge-disjoint) $K_4^{(3)}$ cliques. 
\end{theorem}

From the theory of groups our main tool is a classical result of Erd\H os and Strauss \cite{ES} which states that every finite group contains a large abelian subgroup. The best - and asymptotically optimal - bound is due to Pyber \cite{Py}. 

\begin{theorem}[Pyber]\label{Py}
There is a universal constant $\nu>0$ so that every group of order $n$ contains an abelian subgroup of order at least $e^{\nu\sqrt{\log{n}}}$.
\end{theorem}

Pyber's theorem was also used in a predecessor of this paper, in \cite{So2}. Here we prove a stronger statement which was stated as a conjecture in \cite{So2}. 

\begin{theorem}\label{equ}
For every $\kappa>0$ there is a threshold $n_0\in \mathbb{N}$ such that if $G$ is a finite group of order $|G|\geq n_0$
then the following holds. Any set $H\subset G\times G$ with $|H|\geq \kappa|G|^2$ contains four elements $(\alpha,\beta),$ $(\alpha,\gamma),$ $(\delta,\gamma),$ and $(\delta,\beta)$ such that $\alpha\beta=\delta\gamma.$ 
\end{theorem}

It is easy to see that Theorem \ref{equ} implies Theorem \ref{main}; Every triple $(a,b,ab)$ is uniquely determined by $(a,b)\in   G\times G.$ The triples $(a,b,ab), (a,c,ac),(d,c,dc),$ and $(d,b,db)$ determine at most seven elements of $G$ which are $a,b,c,d, ab=dc, ac,$ and $db.$ (The last two elements might coincide) 

\begin{proof}{(of Theorem \ref{equ})}
Let $A$ be the largest abelian subgroup of $G.$ By Pyber's Theorem we know that $|A|\geq e^{\nu\sqrt{\log{n}}}$. There are elements $\ell, r \in G$ so that $H$ has at least average density in the product of the left and right cosets $\ell A\times Ar,$ that is  $|H\cap (\ell A\times Ar)|\geq \kappa |A|^2.$  

Let us define a 4-partite 3-uniform hypergraph using $H, \ell, r,$ and $A$. The four vertex partitions are $\ell A = V_1, Ar=V_2, \ell Ar=V_3, A=V_4.$  Every triple $(a,b,c)$ where $(a,b)\in H\cap (\ell A\times Ar)$ and $c\in A$ defines four edges, a $K_4^{(3)}$ clique as follows;  $g_i\in V_1, g_j\in V_2, g_k\in V_3,$ and $ g_l\in V_4$ spans a $K_4^{(3)}$ clique if

\begin{itemize}
  \item[1,] $ac=g_i, $
  \item[2,] $cb=g_j,$
  \item[3,] $acb=g_k,$
  \item[4,]  $c=g_l.$
\end{itemize}

With this definition every edge belongs to a unique $(a,b,c)$ triple. As we will see, from the three vertices of an edge one can recover the values of $a, b,$ and $c.$ The number of edges in this hypergraph is at least  $\kappa|A|^3$ and the number of vertices is $4|A|.$ If $|A|$ is large enough in terms of $\kappa$ then by Theorem \ref{FR1} there is a $K_4^{(3)}$ clique in the hypergraph which has edges defined by different 
$(a,b,c)$ triples. Before we continue the proof, let us check that every edge belongs to a unique triple. (For the inverse of an element $g\in G$ we use the usual $g^{-1}$ notation.)

\begin{itemize}
  \item[1,] If  $g_i\in V_1, g_j\in V_2, g_k\in V_3$ spans an edge defined by $(a_1,b_1,c_1)$ then 
  \subitem $a_1=g_kg_j^{-1},$ 
  \subitem $b_1=g_i^{-1}g_k,$ 
  \subitem $c_1=a^{-1}g_i=g_jg_k^{-1}g_i.$ 
  \item[2,] If  $g_i\in V_1, g_j\in V_2, g_l\in V_4$ spans an edge defined by $(a_2,b_2,c_2)$ then 
  \subitem $a_2=g_ig_l^{-1},$ 
  \subitem $b_2=g_l^{-1}g_j,$ 
  \subitem $c_2=g_l.$ 
  \item[3,] If  $g_i\in V_1, g_k\in V_3,g_l\in V_4$ spans an edge defined by $(a_3,b_3,c_3)$ then 
  \subitem $a_3=g_ig_l^{-1},$ 
  \subitem $b_3=g_i^{-1}g_k,$ 
  \subitem $c_3=g_l.$ 
  \item[4,]  If  $ g_j\in V_2, g_k\in V_3,g_l\in V_4$ spans an edge defined by $(a_4,b_4,c_4)$ then 
  \subitem $a_4=g_kg_j^{-1},$ 
  \subitem $b_4=g_l^{-1}g_j,$ 
  \subitem $c_4=g_l.$ 
  \end{itemize}
  
  If two generating triples of the edges of a $K_4^{(3)}$ clique are given, then they determine the vertices and therefore the remaining two edges uniquely. Therefore if a clique is not generated by a single triple then all four edges have distinct generators, $(a_1,b_1,c_1),$  $(a_2,b_2,c_2), $ $(a_3,b_3,c_3),$ and $(a_4,b_4,c_4).$ As we noted before, by Theorem \ref{FR1} we can suppose that such $K_4^{(3)}$ exists if $|A|$ is large enough. 
  
 Note that the four pairs $(a_1,b_1),(a_2,b_2),(a_3,b_3),(a_4,b_4)\in H$ will satisfy the requirements of Theorem \ref{equ}. Set 
 \begin{itemize}
 \item[] $\delta =a_1=a_4,$ 
 \item[] $\beta = b_1=b_3,$
 \item[] $\alpha =a_2=a_3,$ 
 \item[]  $\gamma = b_2=b_4.$
\end{itemize}

It remains to check if $\alpha\beta=\delta\gamma.$ The two elements $c_1$ and $c_2$ are from the abelian subgroup $A,$ so 
\begin{eqnarray*}
c_1^{-1}c_2^{-1}&=&c_2^{-1}c_1^{-1},\\
g_i^{-1}g_kg_j^{-1}g_l^{-1}&=&g_l^{-1}g_i^{-1}g_kg_j^{-1},\\
g_kg_j^{-1}g_l^{-1}g_j&=&g_ig_l^{-1}g_i^{-1}g_k,\\
a_1b_2&=&a_2b_1,\\
\delta\gamma&=&\alpha\beta.
\end{eqnarray*}

\end{proof}

Finally, we briefly bound the number of $(7,4)$-configurations our calculation finds in a group. By the quantitative version of the Frankl-R\"odl theorem, Theroem \ref{FR2}, the number of  $K_4^{(3)}$ cliques for the selected $\ell$ and $r$ elements is at least $c'|A|^4.$ That guarantees at least $c'|A|^3$ $(\alpha,\beta),$ $(\alpha,\gamma),$ $(\delta,\gamma),$ $(\delta,\beta)$ quadruples from $S$ such that $\alpha\beta=\delta\gamma.$ 
Set $S$ has high density in a positive fraction of the left and right cosets.

$$\left |\left \{(r',\ell') : |H\cap (\ell' A\times Ar')|\geq \frac{\kappa}{2} |A|^2\right \}\right |\geq c''\frac{n^2}{|A|^2}.$$  

For these $r',\ell'$ pairs one can repeat the calculations  as we did before, so in each case there are at least $c'''|A|^3$ $(\alpha,\beta),$ $(\alpha,\gamma),$ $(\delta,\gamma),$ $(\delta,\beta)$ quadruples such that $\alpha\beta=\delta\gamma.$ 

\begin{theorem}There is a constant $\mu > 0$ depending on $\kappa,$ the density of $S,$  only so that the number of $(\alpha,\beta),$ $(\alpha,\gamma),$ $(\delta,\gamma),$ $(\delta,\beta)$ quadruples from $S$ such that $\alpha\beta=\delta\gamma$ 
 in the group is at least 
$$\mu |A|n^2\geq \mu  e^{\nu\sqrt{\log{n}}}n^2.$$ 
\end{theorem}
It might be that the right magnitude is $\xi n^3$ for some universal constant $\xi > 0$ independent of the group.

\section{Acknowledgments} 
I am thankful to Vera T. S\'os who continuously encouraged me to work on this problem, and to Noga Alon, Nati Linial, and Endre Szemer\'edi for the useful conversations.

\end{document}